\documentclass[a4paper,10pt,dvipsnames]{article}
\usepackage[utf8]{inputenc}
\usepackage[english]{babel}
\usepackage{amsthm}
\usepackage{amsfonts}
\usepackage{amssymb}	
\usepackage{amsmath}
\allowdisplaybreaks
\usepackage{chemarr}
\usepackage{slashed}
\usepackage{enumerate}
\usepackage{graphicx}
\usepackage{systeme}
\usepackage{dsfont}
\usepackage{relsize}
\usepackage[margin=0.90in]{geometry}
\usepackage{float}
\usepackage{upgreek}
\usepackage{titlesec}
\usepackage{tikz}
\usepackage[labelformat=simple]{subcaption}

\usepackage{graphicx}
\usepackage{bmpsize}
\usepackage{epstopdf}
\usepackage{bbm}
\usepackage{etoolbox}
\usepackage{xcolor}
\usepackage{titling}

\usepackage{blindtext,graphicx}
\usepackage[absolute]{textpos}
\usepackage[hang,flushmargin]{footmisc}

\usepackage{xfrac}
\usepackage{nicefrac}
\usepackage{empheq}
\usepackage{stmaryrd}
\usepackage{cancel}
\usepackage[linguistics]{forest}
\usepackage{url}
\usepackage[font=small]{caption}
\usepackage{array}
\usepackage{dsfont}
\usepackage{tikz}
\usetikzlibrary{trees}
\usetikzlibrary{calc}
\usepackage{pdflscape}
\usepackage{rotating}

\usepackage{todonotes}
\usepackage{xcolor,hyperref}
\definecolor{darkblue}{rgb}{0,0,0.6}
\hypersetup{
    colorlinks=true,       
    linkcolor=darkblue,          
    citecolor=darkblue,        
    filecolor=darkblue,      
    urlcolor=darkblue           
}

\usepackage[giveninits=true,sorting=none,date=year,maxbibnames=99,doi=false,isbn=false,url=false,eprint=false]{biblatex}
\renewbibmacro{in:}{}
\DeclareFieldFormat{pages}{#1}
\AtEveryBibitem{%
\clearlist{language}%
}
\usepackage{csquotes}

\usepackage{setspace}
\makeatletter
\newcommand{\thickhline}{%
\noalign {\ifnum 0=`}\fi \hrule height 1pt

\futurelet \reserved@a \@xhline
}
\newcolumntype{"}{@{\hskip\tabcolsep\vrule width 1pt\hskip\tabcolsep}}
\makeatother
\theoremstyle{definition}

\makeatletter
\newcommand{\pushright}[1]{\ifmeasuring@#1\else\omit\hfill$\displaystyle#1$\fi\ignorespaces}
\newcommand{\pushleft}[1]{\ifmeasuring@#1\else\omit$\displaystyle#1$\hfill\fi\ignorespaces}
\makeatother

\titleformat{\paragraph}
{\normalfont\normalsize\bfseries}{\theparagraph}{1em}{}
\titlespacing*{\paragraph}
{0pt}{3.25ex plus 1ex minus .2ex}{1.5ex plus .2ex}

\newcounter{alphasect}
\def\alphainsection{0}

\let\oldsection=\section
\def\section{%
  \ifnum\alphainsection=1%
    \addtocounter{alphasect}{1}
  \fi%
\oldsection}%
\renewcommand\thesection{%
  \ifnum\alphainsection=1%
    \Alph{alphasect}%
  \else%
    \arabic{section}%
  \fi%
}%

\title{Analytical solutions for long cylindrical shells under radial deformations based on the isotropic relaxed micromorphic continuum}
\author{
Esmaeal Ghavanloo\,\thanks{Corresponding author: Esmaeal Ghavanloo, School of Mechanical Engineering, Shiraz University, Shiraz, 71963-16548, Iran, email: ghavanloo@shirazu.ac.ir},
\quad Pierre Fritsch\,\thanks{Pierre Fritsch, DER de mathématiques, ENS Paris-Saclay, F-91190 Gif-sur-Yvette, France, email: pierre.fritsch@ens-paris-saclay.fr},  \\
and \\
Patrizio Neff\,\thanks{Patrizio Neff, Head of Chair for Nonlinear Analysis and Modelling, Fakultät für Mathematik, Universität Duisburg-Essen, Thea-Leymann-Straße 9, 45127 Essen, Germany, email: patrizio.neff@uni-due.de}
}

\thanksmarkseries{arabic}
\date{\today}

\begin{document}
\maketitle
\begin{abstract}
\noindent
This study presents a closed-form analytical solution for the elastostatic response of long cylindrical shells composed of microstructured materials within the framework of the isotropic relaxed micromorphic continuum. The formulation accounts for microstructural effects by introducing an independent micro-distortion tensor field in addition to the classical displacement field. Under the assumptions of axisymmetric deformation and plane strain conditions, the governing equilibrium equations reduce to a coupled system of ordinary differential equations in the radial coordinate. By introducing suitable auxiliary variables, the system is reformulated into a non-homogeneous modified Bessel equation, which admits an exact analytical solution. Explicit expressions are derived for the radial displacement field and the non-zero components of the micro-distortion tensor. Numerical examples are presented to illustrate the influence of material parameters and the characteristic length on the displacement. The results demonstrate that the relaxed micromorphic model predicts deviations from classical elasticity where microstructural effects are more pronounced. The obtained solution provides valuable physical insight into the mechanics of cylindrical shells and serves as a benchmark for validating numerical implementations of relaxed micromorphic models.
\end{abstract}

{
%
%
%
%
\section{Introduction}
\label{sec:intro}
In natural and synthetic materials that consist of microstructures, such as composite materials \cite{ramirez2025analytical}, metamaterials \cite{abali2022influence}, and porous media \cite{nguyen2022synthesizing}, microstructural effects emerge, giving rise to phenomena that cannot be explained by classical linear elastic continuum mechanics. To describe the effects of microstructures, various types of generalized continuum mechanics models, including micromorphic theories \cite{polizzotto2018micromorphic}, nonlocal elasticity models \cite{shaat2020review}, and higher-order strain/rotation gradient mechanics \cite{faghidian2021unified}, have been developed. These generalized models enable the consideration of long-range interactions \cite{failla2010solution}, capture micro- and nano-scale phenomena \cite{rafii2016nonlocal}, predict material dispersion properties \cite{zeng2006determining}, and facilitate the investigation of realistic material microstructures \cite{arslan2008finite}.

One of the most well-known generalized continuum theories for modeling microstructures is micromorphic continuum mechanics, which was first developed by Eringen \cite{eringen2012microcontinuum} and Mindlin \cite{mindlin1964micro}. In the context of this theory, the deformation consists of macroscopic and microscopic parts. Accordingly, each particle may have independent degrees of freedom describing micro-distortion—nine microscopic degrees of freedom in addition to three macroscopic degrees of freedom \cite{neff2015relaxed}. However, the large number of material parameters (e.g., 18 constants for linear isotropic materials) has so far limited the application of micromorphic continuum mechanics in realistically describing material behavior \cite{madeo2015band}.

To address this limitation, a number of simplified micromorphic formulations have been proposed over the past decade. In 2014, Neff et al. \cite{neff2014unifying} introduced a simplified version of the classical micromorphic theory, known as the relaxed micromorphic model. This formulation substantially reduces the number of material parameters compared to the classical model. Romano et al. \cite{romano2016micromorphic} proposed a non-redundant kinematic description for micromorphic continua by removing the micro-curvature field. They acknowledged that the resulting micromorphic model is non-redundant, in contrast to full micromorphic and linear Cosserat continua. In turn, the relaxed micromorphic model, including the curvature term and assuming a vanishing rotational coupling tensor, constitutes a non-redundant micromorphic formulation. In this setting, the force stress tensor remains symmetric. In contrast, the standard Mindlin–Eringen model and the linear Cosserat model remain redundant. The well-posedness of the relaxed micromorphic model for the case of zero rotational coupling terms is also established using novel generalized Korn’s inequalities for incompatible tensor fields \cite{lewintan2022lp}.  

In 2018, Shaat \cite{shaat2018reduced} introduced the reduced micromorphic model, which was the micro-strain continuum originally developed by Forest and Sievert \cite{forest2006nonlinear}. In 2020, Xiu and Chu \cite{xiu2020micromorphic} developed a modified micromorphic model based on micromechanical principles to describe granular materials. Zhang et al. \cite{zhang2021non} proposed a non-classical Euler–Bernoulli beam model derived from simplified micromorphic elasticity. Their model incorporates microstructural effects to improve predictions of beam behavior at small scales. Recently, Massing et al. \cite{massing2023micromorphic} investigated the applicability of micromorphic models to blood flow in microcirculation, modeling blood as a micromorphic fluid with red blood cells represented as deformable microstructural particles. This approach successfully captures the non-Newtonian behavior of blood, including shear thinning. More recently, the reduced relaxed micromorphic model has been developed and applied to finite-size microstructured metamaterials to describe wave propagation and boundary effects \cite{demetriou2024reduced}.

The mentioned studies have revealed both the versatility of micromorphic-type models and the continuing effort to balance physical realism with analytical and computational simplicity. Among the available formulations, the relaxed micromorphic continuum has emerged as a particularly useful framework for the analysis of microstructured materials. It preserves the essential kinematic features of the classical micromorphic theory while substantially reducing the number of material parameters, resulting in a more manageable yet still physically meaningful model. A key aspect of this formulation is the replacement of the full gradient of the micro-distortion by its Curl without introducing unnecessary constitutive complexity \cite{sky2022primal}. As a result, the relaxed micromorphic model has been successfully employed to analyze a wide range of static \cite{sarhil2023size,rizzi2023analytical,knees2024global} and dynamic problems \cite{barbagallo2019relaxed,d2020effective,rizzi2022boundary} in materials with microstructures. Nevertheless, despite the increasing literature on relaxed micromorphic models, analytical solutions formulated in polar coordinates, especially for axisymmetric problems, remain scarce. In our previous work, we addressed this gap by considering the isotropic relaxed micromorphic model in polar coordinates and deriving explicit analytical solutions for an elastostatic axisymmetric extension problem using modified Bessel functions \cite{ghavanloo2025isotropic}.

In continuation of our previous work, the present study develops analytical solutions for long circular cylindrical shells subjected to radial deformations within the isotropic relaxed micromorphic framework. Cylindrical shells constitute an important class of structures frequently encountered in engineering and applied physics. In such systems, microstructural effects can play a significant role in determining the mechanical response \cite{gao2007variational,collin2009analytical,enakoutsa2018micromorphic}. The problem addressed in this work concerns the elastostatic response of long cylindrical shells made of microstructured materials, modeled within the isotropic relaxed micromorphic continuum, which provides an enriched kinematic description. Under the assumption of axisymmetric deformations, the governing field equations reduce to a coupled system of ordinary differential equations in the radial coordinate. Closed-form analytical solutions are derived using modified Bessel functions. Finally, numerical results are presented to illustrate the influence of key parameters, including the material coefficients of the relaxed micromorphic model and the characteristic length.

%
\section{Brief review of the isotropic relaxed micromorphic model}
According to the general micromorphic theory, the kinematics at each material point are characterized by a macroscopic displacement field ${u} : \Omega \subset \mathbb{R}^3 \rightarrow \mathbb{R}^3$ together with an independent, generally non-symmetric micro-distortion tensor field ${P} : \Omega \subset \mathbb{R}^3 \rightarrow \mathbb{R}^{3 \times 3}$. This formulation enables explicit representation of microstructural deformation mechanisms that classical elasticity cannot capture. The strain energy expression for the isotropic relaxed micromorphic continuum is \cite{gourgiotis2024green}:
{\small
\begin{equation}
\begin{split}
W \left(\text{D}{u}, {P},\mbox{Curl}\,{P}\right) = &
  \, \mu _{\text{e}} \left\lVert \mbox{sym} \left(\text{D}{u} - {P} \right) \right\rVert^{2}
+ \dfrac{\lambda_{\text{e}}}{2} \mbox{tr}^2 \left(\text{D}{u} - {P} \right) 
+ \mu _{\text{c}} \left\lVert \mbox{skew} \left(\text{D}{u} - {P} \right) \right\rVert^{2} \\[3mm]
&
+ \mu_{\text{micro}} \left\lVert \mbox{sym}\,{P} \right\rVert^{2}
+ \dfrac{\lambda_{\text{micro}}}{2} \mbox{tr}^2 \left({P} \right)
+ \dfrac{\mu_{\text{macro}} \,L_{\text{c}}^2 }{2} \, \left\lVert \mbox{Curl} \, {P}\right\rVert^2,
\end{split}
\label{eq:energy}
\end{equation}
}where ($\mu _{\text{e}}$, $\lambda_{\text{e}}$), ($\mu_{\mbox{\tiny micro}}$, $\lambda_{\mbox{\tiny micro}}$), $\mu_c$, $\mu_{\mbox{\tiny macro}}$ and $L_c$ represent, respectively, the meso-scale parameters, the micro-scale parameters, the Cosserat couple modulus,  macro modulus and the characteristic length. In the absence of body forces, the equilibrium equations are \cite{rizzi2021analytical}
{\small
\begin{align}
\text{Div} \, \sigma  =0\, ,
\qquad\qquad\qquad\qquad
\sigma - \sigma_{\rm micro}- \text{Curl}\,m=0\, ,
\label{eq:equi_RM}
\end{align}
}where
{\small
\begin{align}
\label{constEq}
\sigma & \coloneqq
2\mu_{\rm e}\,\text{sym}  (\text{D} u - P )
+ 2\mu_{\rm c}\,\text{skew} (\text{D} u - P )
+ \lambda_{\rm e} \text{tr} (\text{D} u - P ) \mathbbm{1}
\, ,
\notag
\\*
\sigma_{\rm micro} & \coloneqq
2 \mu_{\rm micro}\,\text{sym}\,P
+ \lambda_{\rm micro} \text{tr} \left( P \right) \mathbbm{1}
\, ,
\\*
\notag
m & \coloneqq
\mu_{\rm macro}\, L_{\rm c}^2 \text{Curl} \, P 
\, .
\end{align}}

Here, $\sigma$ denotes the non-symmetric elastic force stress tensor (which becomes symmetric if $\mu_{\rm c}=0$), $\sigma_{\rm micro}$ represents the symmetric microscopic stress tensor, and $m$ corresponds to the non-symmetric moment tensor. Additionally, $\mathbbm{1}$ is the identity matrix, $\text{Div}$ denotes the row-wise divergence operator, $\text{D}$ is the gradient operator, and $\text{sym}$ and $\text{skew}$ refer to the symmetric and skew-symmetric components of a tensor, respectively. Furthermore, Dirichlet and Neumann boundary conditions are 
{\small
\begin{align}
&
\text{Dirichlet:}
&&&&
u
=
\overline{u}
\, ,
&&
\text{and}
&&
P \cdot \tau=\overline{Q} 
\, ,
&
\label{eq:BC_RM_Dir}
\\*
&
\text{Neumann:}
&&&&
t
\coloneqq
\sigma \cdot \, n
= 
0
\, ,
&&
\text{and}
&&
\eta
\coloneqq
m \times n = 
0
\, ,
&
\label{eq:BC_RM_Neu}
\end{align}}
where $\tau$ and $n$ denote the outward tangent and normal vectors on the boundary, respectively. The higher-order Dirichlet boundary conditions in Eq. (\ref{eq:BC_RM_Dir}) can be specified as follows \cite{schroder2022lagrange}: 
{\small
\begin{align}
P \cdot \tau
=
\overline{Q}
=
\text{D}u \cdot \tau
\, ,
\label{eq:CCBC_RM}
\end{align}}
referred to as consistent coupling boundary conditions \cite{d2022consistent}. 

\section{Deformation of long cylindrical shells}
This section addresses the axisymmetric deformation of a long, circular cylindrical shell with inner radius $r_i$ and outer radius $r_o$ (Figure 1). The displacement field and the micro-distortion tensor are functions solely of the radial coordinate, represented by $u_r(r)$, $P_{rr}(r)$, $P_{\theta\theta}(r)$, $P_{r\theta}(r)$, and $P_{\theta r}(r)$, with the circumferential displacement component $u_\theta$ identically zero. To enhance readability, we adopt the following shorthand notation: $\mu_{\rm M} \equiv \mu_{\rm macro}$, $\mu_{\rm m} \equiv \mu_{\rm micro}$, $\lambda_{\rm m} \equiv \lambda_{\rm micro}$, and $\lambda_{\rm M} \equiv \lambda_{\rm macro}$. The governing equations are given by \cite{ghavanloo2025isotropic}:
{\small
\begin{align}
  (2{{\mu }_{\rm e}}+{{\lambda }_{\rm e}})\frac{d}{dr}\left( \frac{d{{u}_{r}}}{dr}-{{P}_{rr}} \right)+{{\lambda }_{\rm e}}\frac{d}{dr}\left( \frac{{{u}_{r}}}{r}-{{P}_{\theta \theta }} \right)=-2{{\mu }_{\rm e}}\frac{d}{dr}\left( \frac{{{u}_{r}}}{r}-{{P}_{\theta \theta }} \right)-2{{\mu }_{\rm e}}\left( \frac{d{{P}_{\theta \theta }}}{dr}+\frac{{{P}_{\theta \theta }}-{{P}_{rr}}}{r} \right) ,
  \label{EAEP1}
\end{align}

\begin{align}
(2{{\mu }_{\rm e}}+{{\lambda }_{\rm e}})\left( \frac{d{{u}_{r}}}{dr}-{{P}_{rr}} \right)+{{\lambda }_{\rm e}}\left( \frac{{{u}_{r}}}{r}-{{P}_{\theta \theta }} \right)=2{{\mu }_{\rm m}}{{P}_{rr}}+{{\lambda }_{\rm m}}({{P}_{rr}}+{{P}_{\theta \theta }})-\frac{{{\mu }_{\textrm{M}}}L_{c}^{2}}{r}\left( \frac{d{{P}_{\theta \theta }}}{dr}+\frac{{{P}_{\theta \theta }}-{{P}_{rr}}}{r} \right) ,
\label{EAEP3}
\end{align}

\begin{align}
(2{{\mu }_{\rm e}}+{{\lambda }_{\rm e}})\left( \frac{{{u}_{r}}}{r}-{{P}_{\theta \theta }} \right)+{{\lambda }_{\rm e}}\left( \frac{d{{u}_{r}}}{dr}-{{P}_{rr}} \right)=2{{\mu }_{\rm m}}{{P}_{\theta \theta }}+{{\lambda }_{\rm m}}({{P}_{rr}}+{{P}_{\theta \theta }})-{{\mu }_{\textrm{M}}}L_{c}^{2}\frac{d}{dr}\left( \frac{d{{P}_{\theta \theta }}}{dr}+\frac{{{P}_{\theta \theta }}-{{P}_{rr}}}{r} \right)  .
\label{EAEP6}
\end{align}

\begin{align}
{{\mu }_{\text{e}}}\left( \frac{d{{P}_{r\theta }}}{dr}+\frac{d{{P}_{\theta r}}}{dr}+2\frac{{{P}_{r\theta }}+{{P}_{\theta r}}}{r} \right)={{\mu }_{\text{c}}}\left( \frac{d{{P}_{r\theta }}}{dr}-\frac{d{{P}_{\theta r}}}{dr} \right),
  \label{EAEP2}
\end{align}

\begin{align}
{{\mu }_{\text{e}}}\left( {{P}_{r\theta }}+{{P}_{\theta r}} \right)+{{\mu }_{\text{m}}}({{P}_{r\theta }}+{{P}_{\theta r}})+{{\mu }_{\text{c}}}\left( {{P}_{r\theta }}-{{P}_{\theta r}} \right)+{{\mu }_{\text{M}}}L_{\text{c}}^{\text{2}}\left( \frac{1}{r}\frac{d{{P}_{\theta r}}}{dr}+\frac{{{P}_{r\theta }}+{{P}_{\theta r}}}{{{r}^{2}}} \right)=0,
\label{EAEP4}
\end{align}

\begin{align}
{{\mu }_{\text{e}}}\left( {{P}_{r\theta }}+{{P}_{\theta r}} \right)+{{\mu }_{\text{m}}}({{P}_{r\theta }}+{{P}_{\theta r}})+{{\mu }_{\text{c}}}\left( {{P}_{\theta r}}-{{P}_{r\theta }} \right)+{{\mu }_{\text{M}}}L_{\text{c}}^{\text{2}}\left( -\frac{{{d}^{2}}{{P}_{\theta r}}}{d{{r}^{2}}}+\frac{{{P}_{r\theta }}+{{P}_{\theta r}}}{{{r}^{2}}}-\frac{1}{r}\frac{d{{P}_{r\theta }}}{dr}-\frac{1}{r}\frac{d{{P}_{\theta r}}}{dr} \right)=0,
\label{EAEP5}
\end{align}}

\begin{figure}[!ht]
\centering
\includegraphics[scale=0.2]{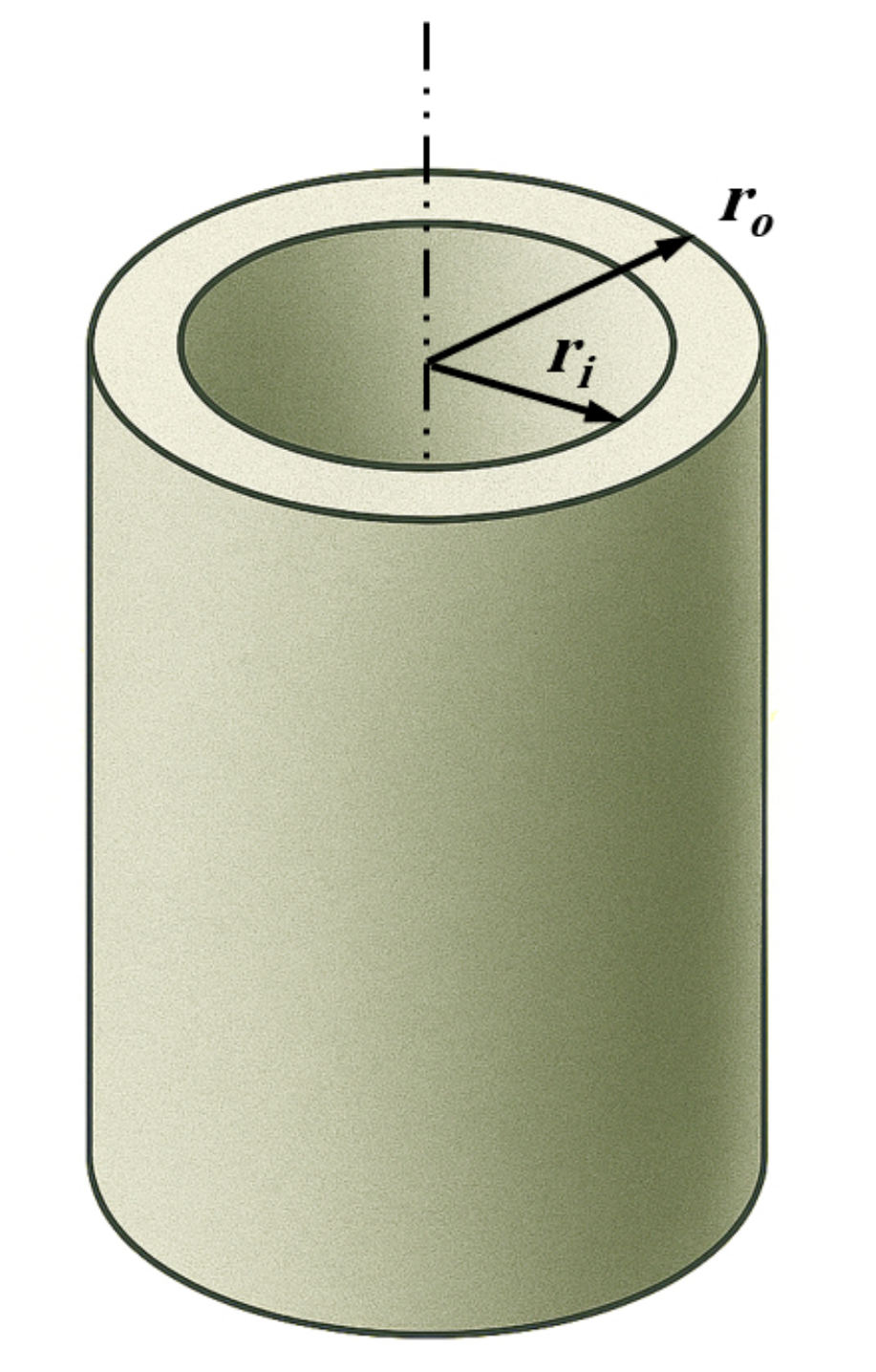}
\caption{A long cylindrical shell.}
\label{figure1}
\end{figure}

The Dirichlet boundary conditions for the problem are
{\small
\begin{align}
{{u}_{r}}(r_i)={{U}_{i}}.
\label{boudcon1}
\end{align}

\begin{align}
{{u}_{r}}(r_o)={{U}_{o}}.
\label{boudcon2}
\end{align}}

The choice of boundary conditions is motivated by the objective of deriving a benchmark analytical solution to evaluate the robustness of a numerical implementation of the relaxed micromorphic model within a finite element framework. This task is nontrivial because the model involves higher-order kinematics, which are not commonly supported in standard finite element codes \cite{sky2022primal,sarhil2023size}. 

Furthermore, the consistent coupling boundary conditions are as follows \cite{ghavanloo2025isotropic}: 
{\small
\begin{align}
{{P}_{\theta \theta }}(r_i)=\frac{{{U}_{i}}}{r_i},\quad \quad {{P}_{\theta \theta }}(r_o)=\frac{{{U}_{o}}}{r_o},
\label{boudcon3}
\end{align}}
{\small\begin{align}
{{P}_{r\theta }}(r_i)={{P}_{r\theta }}(r_o)=0. 
\label{boudcon4}
\end{align}}}

It should be noted that, although there are five unknown variables and six equations, the governing equations naturally split into two independent subsystems due to axisymmetry. The first subsystem involves the radial quantities 
$u_r$, $P_{rr}$, and $P_{\theta\theta}$, while the second subsystem governs $P_{r\theta}$ and $P_{\theta r}$. The latter subsystem is homogeneous and decoupled from the radial fields. As will be shown later, under the consistent coupling boundary conditions and axisymmetric kinematics, this subsystem admits only the trivial solution $P_{r\theta} = P_{\theta r} = 0.$ Consequently, the effective problem reduces to the determination of the radial displacement and the diagonal components of the micro-distortion tensor.

Because the shell is long, plane-strain conditions are justified, and the bulk micro-moduli $\kappa_{\rm e}$ and $\kappa_{\rm m} \equiv \kappa_{\rm micro}$ relate to their respective Lamé-type micro-moduli through the following two-dimensional expressions \cite{gourgiotis2024green}:
{\small
\begin{equation}\label{eq:KAPPA}
\kappa_{\rm e}\coloneqq \lambda_{\rm e}+\mu_{\rm e}, \qquad \quad \kappa_{\rm m} \coloneqq\lambda_{\rm m}+\mu_{\rm m}\,.
\end{equation}}

Furthermore, the relationships between the macro moduli and the corresponding micro-moduli can be expressed using the \textit{Reuss-like homogenization formula} as follows \cite{rizzi2022analytical}:
{\small
\begin{align}
\label{defmod}
\mu_{\rm M} &\coloneqq \dfrac{\mu_{\rm e} \, \mu_{\rm m}}{\mu_{\rm e} + \mu_{\rm m}} 
\qquad\Leftrightarrow\qquad \frac{1}{\mu_{\rm M}}=\frac{1}{\mu_{\rm e}}+\frac{1}{\mu_{\rm m}}\,,\notag\\
\kappa_{\rm M} &\coloneqq \dfrac{\kappa_{\rm e} \, \kappa_{\rm m}}{\kappa_{\rm e} + \kappa_{\rm m}} \qquad\Leftrightarrow\qquad
 \frac{1}{\kappa_{\rm M}}=\frac{1}{\kappa_{\rm e}}+\frac{1}{\kappa_{\rm m}}\,,
 \end{align}
}
where $\kappa_{\rm M}= \mu_{\rm M}+\lambda_{\rm M} $. It should be noted that $\kappa_{\rm M}$ and $\mu_{\rm M}$ can be uniquely determined using classical homogenization theory for a given periodic microstructured material. This feature distinguishes the relaxed micromorphic model from other generalized continuum models. The requirement of positive definiteness of the strain energy density imposes the following constraints on the material parameters \cite{neff2014unifying}:

\begin{align}
\mu_{\rm m} &> 0, \qquad
\mu_{\rm e} > 0, \qquad
\mu_{\rm c} \geq 0, \qquad
\kappa_{\rm m} > 0, \qquad
\kappa_{\rm e} > 0.
\label{constraints}
\end{align}

Using equation (\ref{defmod}), these positivity conditions imply the inequalities $\mu_{\rm e} > \mu_{\rm M} > 0$ and $\kappa_{\rm e} > \kappa_{\rm M} > 0$. Moreover, for physically admissible isotropic materials under plane–strain conditions, we have $\kappa_{\rm M} > \mu_{\rm M}$ and $\kappa_{\rm m} > \mu_{\rm m}$.

Deriving the analytical solution requires a reformulation of the system given by equations (\ref{EAEP1})–(\ref{EAEP6}). As a first step, equation (\ref{EAEP1}) is rearranged in the following form:
{\small
\begin{align}
  & ({{\kappa }_{\rm e}}+{{\mu }_{\rm e}})\frac{d}{dr}\left( \frac{d{{u}_{r}}}{dr}+\frac{{{u}_{r}}}{r} \right)-2{{\mu }_{\rm e}}(\frac{d{{P}_{rr}}}{dr}-\frac{{{P}_{\theta \theta }}-{{P}_{rr}}}{r})-({{\kappa }_{\rm e}}-{{\mu }_{\rm e}})\frac{d}{dr}({{P}_{\theta \theta }}+{{P}_{rr}})=0 .
  \label{EAEP_ref}
\end{align}
}

Next, differentiating equation (\ref{EAEP3}) with respect to $r$, and substituting equations (\ref{EAEP1}) and (\ref{EAEP6}) along with the applied simplifications, yields the following result:
{\small
\begin{align}
2{{\mu }_{\rm m}}\left( \frac{d{{P}_{rr}}}{dr}-\frac{{{P}_{\theta \theta }}-{{P}_{rr}}}{r} \right)+({{\kappa }_{\rm m}}-{{\mu }_{\rm m}})\frac{d}{dr}({{P}_{\theta \theta }}+{{P}_{rr}})=0 .
\label{ref_EAEP2}
\end{align}}

Substitution of equation (\ref{ref_EAEP2}) into equation (\ref{EAEP_ref}) gives:

{\small
\begin{align}
\frac{d{{u}_{r}}}{dr}+\frac{{{u}_{r}}}{r}=\frac{{{\mu }_{\rm m}}{{\kappa }_{\rm e}}-{{\mu }_{\rm e}}{{\kappa }_{\rm m}}}{{{\mu }_{\rm m}}({{\mu }_{\rm e}}+{{\kappa}_{\rm e}})}({{P}_{\theta \theta }}+{{P}_{rr}})+{{C}_{1}} .
\label{Int_EAEP}
\end{align}}

Here, \( C_1 \) denotes an integration constant. Furthermore, adding equations (\ref{EAEP3}) and (\ref{EAEP6}) yields:
{\small
\begin{align}
2{{\kappa }_{\rm e}}\left( \frac{d{{u}_{r}}}{dr}+\frac{{{u}_{r}}}{r} \right)=2({\kappa}_{\rm m}+{\kappa}_{\rm e})({{P}_{rr}}+{{P}_{\theta \theta }})-{{\mu }_{\rm M}}L_{c}^{2}\frac{d}{dr}\left( \frac{d{{P}_{\theta \theta }}}{dr}+\frac{{{P}_{\theta \theta }}-{{P}_{rr}}}{r} \right)-{{\mu }_{\rm M}}L_{c}^{2}\frac{1}{r}\left( \frac{d{{P}_{\theta \theta }}}{dr}+\frac{{{P}_{\theta \theta }}-{{P}_{rr}}}{r} \right).
\label{EAEP3_ref}
\end{align}}

At this stage, the following variables are introduced to facilitate the analysis:
{\small
\begin{align}
X=\frac{d{{u}_{r}}}{dr}-{{P}_{rr}}, \quad \quad \quad Y=\frac{{{u}_{r}}}{r}-{{P}_{\theta \theta }}, \quad \quad \quad Z={{P}_{\theta \theta }}+{{P}_{rr}} .
\label{new_var}
\end{align}}

By substituting equation (\ref{new_var}) into equations (\ref{ref_EAEP2})–(\ref{EAEP3_ref}), the following expressions are obtained:
{\small
\begin{align}
\frac{dY}{dr}+\frac{Y-X}{r}&=-\left(\frac{{{\kappa}_{\rm m}}+{{\mu }_{\rm m}}}{2{{\mu }_{\rm m}}}\right)\frac{dZ}{dr}, \\
X+Y&=-\frac{{{\mu }_{\rm e}}({{\kappa}_{\rm m}}+{{\mu }_{\rm m}})}{{{\mu }_{\rm m}}({{\kappa}_{\rm e}}+{{\mu }_{\rm e}})}Z+{{C}_{1}},\\
X+Y&=\frac{{{\kappa}_{\rm m}}}{{{\kappa}_{\rm e}}}Z+\frac{{{\mu }_{\rm M}}L_{c}^{2}}{2{{\kappa}_{\rm e}}}\frac{d}{dr}\left( \frac{dY}{dr}+\frac{Y-X}{r}\right)+\frac{{{\mu }_{\rm M}}L_{c}^{2}}{2{{\kappa}_{\rm e}}}\frac{1}{r}\left( \frac{dY}{dr}+\frac{Y-X}{r} \right).
\end{align}}

Using equations (25) and (26), and after simplification, equation (27) can be rewritten as:
{\small
\begin{align}
\frac{{{d}^{2}}Z}{d{{r}^{2}}}+\frac{1}{r}\frac{dZ}{dr}-aZ+b=0,
\label{EQZ}
\end{align}}
where
{\small
\begin{align}
a=\frac{4}{{{\mu }_{\rm M}}L_{c}^{2}}\left( \frac{{{\mu }_{\rm e}}{{\kappa}_{\rm e}}}{{{\kappa}_{\rm e}}+{{\mu}_{\rm e}}}+\frac{{{\mu }_{\rm m}}{{\kappa}_{\rm m}}}{{{\kappa}_{\rm m}}+{{\mu}_{\rm m}}} \right), \quad \quad b=\frac{4}{{{\mu }_{\rm M}}L_{c}^{2}}\frac{{{C}_{1}}{{\kappa}_{\rm e}}{{\mu }_{\rm m}}}{{{\kappa}_{\rm m}}+{{\mu }_{\rm m}}}.
\label{equa}
\end{align}}

Equation (\ref{EQZ}) is a non-homogeneous modified Bessel equation, and its general solution is given by:
{\small
\begin{align}
Z(r)=\frac{b}{a}+{{D}_{1}}{{I}_{0}}(\sqrt{a}r)+{{D}_{2}}{{K}_{0}}(\sqrt{a}r),
\label{Sol}
\end{align}}
in which $D_1$ and $D_2$ are constants, while $I_0$ and $K_0$ denote the modified Bessel functions of the first and second kind, respectively, of order zero. By applying equations (\ref{new_var}) and (\ref{Sol}), equation (\ref{Int_EAEP}) can be reformulated as:
{\small
\begin{align}
\frac{d{{u}_{r}}}{dr}+\frac{{{u}_{r}}}{r}={{C}_{1}}A+B({{D}_{1}}{{I}_{0}}(\sqrt{a}r)+{{D}_{2}}{{K}_{0}}(\sqrt{a}r)), 
\label{EQu}
\end{align}}
where
{\small
\begin{align}
A=\frac{{{\mu }_{\rm m}}({{\kappa }_{\rm e}}+{{\mu }_{\rm e}})\left( {{\kappa }_{\rm e}}+{{\kappa }_{\rm m}} \right)}{{{\mu }_{\rm e}}{{\kappa }_{\rm e}}({{\mu }_{\rm m}}+{{\kappa }_{m}})+{{\mu }_{\rm m}}{{\kappa }_{\rm m}}({{\mu }_{\rm e}}+{{\kappa }_{\rm e}})} , \quad \quad \quad B=\frac{{{\mu }_{\rm m}}{{\kappa }_{\rm e}}-{{\mu }_{\rm e}}{{\kappa }_{\rm m}}}{{{\mu }_{\rm m}}({{\mu }_{\rm e}}+{{\kappa }_{\rm e}})}.
\label{A_B}
\end{align}
}

The general solution of equation (\ref{EQu}) can be obtained by:
{\small
\begin{eqnarray}
u_r=\frac{{{C}_{1}}{{\mu }_{\rm m}}r\left( {{\kappa }_{\rm e}}+{{\kappa }_{m}} \right)\left( {{\kappa }_{\rm e}}+{{\mu }_{\rm e}} \right)}{2\left( {{\kappa }_{\rm e}}{{\mu }_{\rm e}}\left( {{\kappa }_{\rm m}}+{{\mu }_{\rm m}} \right)+{{\kappa }_{\rm m}}{{\mu }_{\rm m}}\left( {{\kappa }_{\rm e}}+{{\mu }_{\rm e}} \right) \right)}+\frac{{{C}_{2}}}{r}+\frac{{{\kappa }_{\rm e}}{{\mu }_{\rm m}}-{{\kappa }_{\rm m}}{{\mu }_{\rm e}}}{\sqrt{a}{{\mu }_{\rm m}}\left( {{\kappa }_{\rm e}}+{{\mu }_{\rm e}} \right)}\left( {{D}_{1}}{{I}_{1}}\left( \sqrt{a}r \right)-{{D}_{2}}{{K}_{1}}\left( \sqrt{a}r \right) \right),
\label{Solution_u}
\end{eqnarray}}
in which $C_2$ is an integration constant, and $I_1$ denotes the modified Bessel function of the first kind of order one. 

Closed-form expressions for the non-zero components of the micro-distortion tensor $P$ are derived by summing equations (25) and (26). The resulting equation is:
{\small
\begin{align}
\frac{dY}{dr}+\frac{2Y}{r}=-{{D}_{1}}{{\xi }_{1}}\sqrt{a}{{I}_{1}}(\sqrt{a}r)+\frac{{{D}_{1}}{{\xi }_{2}}{{I}_{0}}(\sqrt{a}r)}{r}+{{D}_{2}}{{\xi }_{1}}\sqrt{a}{{K}_{1}}(\sqrt{a}r)+\frac{{{D}_{2}}{{\xi }_{2}}{{K}_{0}}(\sqrt{a}r)}{r}+\frac{{{C}_{1}}{{\xi }_{3}}}{r},
\label{dY}
\end{align}}
where
{\small
\begin{align}
{{\xi }_{1}}=\frac{{{\kappa }_{\rm m}}+{{\mu }_{\rm m}}}{2{{\mu }_{\rm m}}}, \quad \quad \quad {{\xi }_{2}}=-\frac{{{\mu }_{\rm e}}({{\kappa }_{\rm m}}+{{\mu }_{\rm m}})}{{{\mu }_{\rm m}}({{\kappa }_{\rm e}}+{{\mu }_{\rm e}})}, \quad \quad \quad {{\xi }_{3}}=\frac{{{\mu }_{\rm m}}{{\kappa }_{\rm m}}({{\mu }_{\rm e}}+{{\kappa }_{\rm e}})}{{{\mu }_{\rm e}}{{\kappa }_{\rm e}}({{\mu }_{\rm m}}+{{\kappa }_{\rm m}})+{{\mu }_{\rm m}}{{\kappa }_{\rm m}}({{\mu }_{\rm e}}+{{\kappa }_{\rm e}})}.
\label{zeta}
\end{align}
}

The solution of Eq. (\ref{dY}) is
{\small
\begin{align}
Y=\frac{{{C}_{1}}{{\xi }_{3}}}{2}+\frac{{{C}_{3}}}{{{r}^{2}}}-{{D}_{1}}\left( {{\xi }_{1}}{{I}_{0}}(\sqrt{a}r)-\frac{{{I}_{1}}(\sqrt{a}r)}{\sqrt{a}r}\left( 2{{\xi }_{1}}+{{\xi }_{2}} \right) \right)-{{D}_{2}}\left( {{\xi }_{1}}{{K}_{0}}(\sqrt{a}r)+\frac{{{K}_{1}}(\sqrt{a}r)}{\sqrt{a}r}\left( 2{{\xi }_{1}}+{{\xi }_{2}} \right) \right),
\label{SolY}
\end{align}
{\small
where $C_3$ is an integration constant. Using equations (\ref{new_var}), (\ref{Solution_u}), and (\ref{SolY}), the diagonal components of the micro-distortion tensor can be determined as follows:
{\small
\begin{eqnarray}
{{P}_{\theta \theta }}&=&\frac{{{C}_{1}}{{\mu }_{\rm m}}{{\kappa }_{\rm e}}({{\kappa }_{\rm e}}+{{\mu }_{\rm e}})}{2({{\mu }_{\rm e}}{{\kappa }_{\rm e}}({{\mu }_{\rm m}}+{{\kappa }_{\rm m}})+{{\mu }_{\rm m}}{{\kappa }_{\rm m}}({{\mu }_{\rm e}}+{{\kappa }_{\rm e}}))}+{{D}_{1}}\left( \frac{{{\kappa }_{\rm m}}+{{\mu }_{\rm m}}}{2{{\mu }_{\rm m}}}{{I}_{0}}(\sqrt{a}r)-\frac{{{I}_{1}}(\sqrt{a}r)}{\sqrt{a}r}\frac{{{\kappa }_{\rm m}}}{{{\mu }_{\rm m}}} \right) \\ \nonumber
&&+{{D}_{2}}\left( \frac{{{\kappa }_{\rm m}}+{{\mu }_{\rm m}}}{2{{\mu }_{\rm m}}}{{K}_{0}}(\sqrt{a}r)+\frac{{{K}_{1}}(\sqrt{a}r)}{\sqrt{a}r}\frac{{{\kappa }_{\rm m}}}{{{\mu }_{\rm m}}} \right)+\frac{1}{{{r}^{2}}}\left( {{C}_{2}}-{{C}_{3}} \right),
\label{Ptt}
\end{eqnarray}
\begin{eqnarray}
{{P}_{rr}}&=&\frac{{{C}_{1}}{{\kappa }_{\rm e}}{{\mu }_{\rm m}}({{\mu }_{\rm e}}+{{\kappa }_{\rm e}})}{2({{\mu }_{\rm e}}{{\kappa }_{\rm e}}({{\mu }_{\rm m}}+{{\kappa }_{\rm m}})+{{\mu }_{\rm m}}{{\kappa }_{\rm m}}({{\mu }_{\rm e}}+{{\kappa }_{\rm e}}))}-{{D}_{1}}\left( \frac{{{\kappa }_{\rm m}}-{{\mu }_{\rm m}}}{2{{\mu }_{\rm m}}}{{I}_{0}}(\sqrt{a}r)-\frac{{{I}_{1}}(\sqrt{a}r)}{\sqrt{a}r}\frac{{{\kappa }_{\rm m}}}{{{\mu }_{\rm m}}} \right) \\ \nonumber
&&-{{D}_{2}}\left( \frac{{{\kappa }_{\rm m}}-{{\mu }_{\rm m}}}{2{{\mu }_{\rm m}}}{{K}_{0}}(\sqrt{a}r)+\frac{{{K}_{1}}(\sqrt{a}r)}{\sqrt{a}r}\frac{{{\kappa }_{\rm m}}}{{{\mu }_{\rm m}}} \right)-\frac{1}{{{r}^{2}}}\left( {{C}_{2}}-{{C}_{3}} \right).
\label{Prr}
\end{eqnarray}}

It is important to note that five unknown parameters ($C_1$, $C_2$, $C_3$, $D_1$, and $D_2$) must be determined, while only four boundary conditions are available (\ref{boudcon1}-\ref{boudcon3}).
An additional condition is obtained by substituting equations (\ref{Solution_u}), (\ref{Ptt}), and (\ref{Prr}) into equation (\ref{EAEP3}), resulting in:
{\small
\begin{eqnarray}
{{C}_{3}}={{C}_{2}}\frac{{{\mu }_{\rm m}}}{{{\mu }_{\rm e}}+{{\mu }_{\rm m}}}.
\label{cond4}
\end{eqnarray}}

Now, using boundary conditions (\ref{boudcon1}-\ref{boudcon3}) and Eq. (\ref{cond4}), we obtain
{\small
\begin{eqnarray}
\left[ \begin{matrix}
   {{F}_{0}} & \frac{1}{r_{i}^{2}} & {{F}_{1}}({{r}_{i}}) & {{F}_{2}}({{r}_{i}})  \\
   {{F}_{0}} & \frac{1}{r_{o}^{2}} & {{F}_{1}}({{r}_{o}}) & {{F}_{2}}({{r}_{o}})  \\
   \frac{{{F}_{0}}{{\kappa }_{e}}}{\left( {{\kappa }_{e}}+{{\kappa }_{m}} \right)} & \frac{1}{r_{i}^{2}}\left( \frac{{{\mu }_{e}}}{{{\mu }_{e}}+{{\mu }_{m}}} \right) & {{F}_{3}}({{r}_{i}}) & {{F}_{4}}({{r}_{i}})  \\
   \frac{{{F}_{0}}{{\kappa }_{e}}}{\left( {{\kappa }_{e}}+{{\kappa }_{m}} \right)} & \frac{1}{r_{o}^{2}}\left( \frac{{{\mu }_{e}}}{{{\mu }_{e}}+{{\mu }_{m}}} \right) & {{F}_{3}}({{r}_{o}}) & {{F}_{4}}({{r}_{o}})  \\
\end{matrix} \right]\left( \begin{matrix}
   {{C}_{1}}  \\
   {{C}_{2}}  \\
   {{D}_{1}}  \\
   {{D}_{2}}  \\
\end{matrix} \right)=\left( \begin{matrix}
   \frac{{{U}_{i}}}{{{r}_{i}}}  \\
   \frac{{{U}_{o}}}{{{r}_{o}}}  \\
   \frac{{{U}_{i}}}{{{r}_{i}}}  \\
   \frac{{{U}_{o}}}{{{r}_{o}}}  \\
\end{matrix} \right)
\label{Conditionmat}
\end{eqnarray}}
where
{\small
\begin{eqnarray}
{{F}_{0}}&=&\frac{{{\mu }_{\rm m}}\left( {{\kappa }_{\rm e}}+{{\kappa }_{\rm m}} \right)\left( {{\kappa }_{\rm e}}+{{\mu }_{\rm e}} \right)}{2\left( {{\kappa }_{\rm e}}{{\mu }_{\rm e}}\left( {{\kappa }_{\rm m}}+{{\mu }_{\rm m}} \right)+{{\kappa }_{\rm m}}{{\mu }_{\rm m}}\left( {{\kappa }_{\rm e}}+{{\mu }_{\rm e}} \right) \right)}, \\ \nonumber
{{F}_{1}}(r)&=&\frac{{{\kappa }_{\rm e}}{{\mu }_{\rm m}}-{{\kappa }_{\rm m}}{{\mu }_{\rm e}}}{\sqrt{a}{{\mu }_{\rm m}}\left( {{\kappa }_{\rm e}}+{{\mu }_{\rm e}} \right)}{{I}_{1}}\left( \sqrt{a}r \right), \\ \nonumber
{{F}_{2}}(r)&=&-\frac{{{\kappa }_{\rm e}}{{\mu }_{\rm m}}-{{\kappa }_{\rm m}}{{\mu }_{\rm e}}}{\sqrt{a}{{\mu }_{\rm m}}\left( {{\kappa }_{\rm e}}+{{\mu }_{\rm e}} \right)r}{{K}_{1}}\left( \sqrt{a}r \right), \\ \nonumber
{{F}_{3}}(r)&=&\frac{{{\kappa }_{\rm m}}+{{\mu }_{\rm m}}}{2{{\mu }_{\rm m}}}{{I}_{0}}(\sqrt{a}r)-\frac{{{I}_{1}}(\sqrt{a}r)}{\sqrt{a}r}\frac{{{\kappa }_{\rm m}}}{{{\mu }_{\rm m}}}, \\ \nonumber
{{F}_{4}}(r)&=&\frac{{{\kappa }_{\rm m}}+{{\mu }_{\rm m}}}{2{{\mu }_{\rm m}}}{{K}_{0}}(\sqrt{a}r)+\frac{{{K}_{1}}(\sqrt{a}r)}{\sqrt{a}r}\frac{{{\kappa }_{\rm m}}}{{{\mu }_{\rm m}}}.
\label{F}
\end{eqnarray}}

Solving the linear system given in equation (\ref{Conditionmat}) yields the values of the remaining unknown coefficients $C_1$, $C_2$, $D_1$, and $D_2$. These coefficients fully define the displacement and the diagonal components of the micro-distortion tensor, subject to the boundary conditions and the additional constraint provided by equation (\ref{cond4}). 

To complete the solution procedure, it is essential to derive relations for $P_{r\theta}$ and $P_{\theta r}$. In this context, subtracting equation (\ref{EAEP4}) from equation (\ref{EAEP5}) yields the following result:
{\small
\begin{align}
2{{\mu }_{\text{c}}}\left( {{P}_{r\theta }}-{{P}_{\theta r}} \right)=-\frac{{{\mu }_{\text{M}}}L_{\text{c}}^{2}}{r}\left( \frac{d}{dr}\left( r\frac{d{{P}_{\theta r}}}{dr}+{{P}_{\theta r}}+{{P}_{r\theta }} \right) \right).
\label{SEq65}
\end{align}
}
Using equations (\ref{EAEP2}) and (\ref{EAEP5}), equation (\ref{SEq65}) can be rearranged in the following form:
{\small
\begin{align}
\frac{d}{dr}\left( {{P}_{r\theta }} + {{P}_{\theta r}} \right) + \frac{2}{r} \left( {{P}_{r\theta }} + {{P}_{\theta r}} \right) = 0.
\label{newp}
\end{align}}
The solution of equation (\ref{newp}) is
{\small
\begin{align}
{{P}_{r\theta }}+{{P}_{\theta r}}=\frac{{{C}_{4}}}{{{r}^{2}}}.
\label{solnewp}
\end{align}}

Substituting equation (\ref{newp}) into Eq. (\ref{EAEP2}) yields:
\begin{align}
{{P}_{r\theta }}-{{P}_{\theta r}}={{C}_{5}},
\label{solnewp2}
\end{align}}
where $C_4$ and $C_5$ are integration constants. By applying the consistent coupling boundary conditions (\ref{boudcon4}), we find ${{P}_{r\theta }}={{P}_{\theta r}}=0$. Finally, it is important to note that in the axisymmetric case, all terms involving the Cosserat couple modulus $\mu_{\rm c}$ vanish. Consequently, $\mu_{\rm c}$ does not appear in the final solutions.

\section{Numerical results and discussion}
In this section, numerical examples are provided to demonstrate the analytical solution obtained in the previous section and to examine the effects of material parameters, boundary displacements, and geometric properties on the mechanical response of cylindrical shells within the relaxed micromorphic model. The mechanical behavior can be fully described using four material parameters. To provide a unified framework for all the following numerical results, the dimensionless quantities 
$G_i$ ($i$=1,2,3) are introduced:

\begin{eqnarray}
\mu_{\rm m} = G_1 \,\mu_{\rm M}, \quad \quad \kappa_{\rm m} = G_2 \,\mu_{\rm M}, \quad \quad \kappa_{\rm M} = G_3 \,\mu_{\rm M},
\end{eqnarray}
where, based on the inequalities satisfied by the material parameters, we have $G_2 > G_1> 1$ and $G_2 > G_3> 1$.

To investigate the influence of the shear modulus ratio $G_1$, numerical results are presented for two representative geometric configurations: a thick cylindrical shell and a thin cylindrical shell. This distinction allows examination of the sensitivity of the mechanical response to microstructural effects across different thickness regimes. Figure~\ref{figure2} illustrates the normalized radial displacement for a thick shell with an inner-to-outer radius ratio $\beta=0.15$, $G_1 = 1.45, 3.25, 4.95$, $r_o/L_{\rm c}=2$, $G_2=5$, $G_3=2$, and $U_i/U_o=0$. Similarly, Figure~\ref{figure3} presents the normalized radial displacement for a thin shell with the same material and loading parameters as in Figure~\ref{figure2}, except that the inner-to-outer radius ratio is $\beta=0.85$. For the thick shell, when $G_1 = 1.45$, the relaxed micromorphic model predicts normalized radial displacements that exceed those of classical linear elasticity. This indicates that, at lower micro shear moduli, the additional degrees of freedom inherent to the micromorphic framework produce a more compliant structural response compared to the classical prediction. In contrast, increasing $G_1$ to $3.25$ and $4.95$ reverses this trend, with predicted displacements falling below the classical solution. This behavior reflects a stiffening effect, whereby higher micro shear moduli increase resistance to deformation, which is not captured by classical theory. While these deviations are pronounced in the thick shell ($r/R_o = 0.15$), the thin shell ($r/R_o = 0.85$) is largely insensitive to variations in $G_1$, with all results converging toward the classical linear elasticity solution due to geometric constraints. Accordingly, the following numerical results focus exclusively on the thick shell.

\begin{figure}[!ht]
\centering
\includegraphics[scale=0.3]{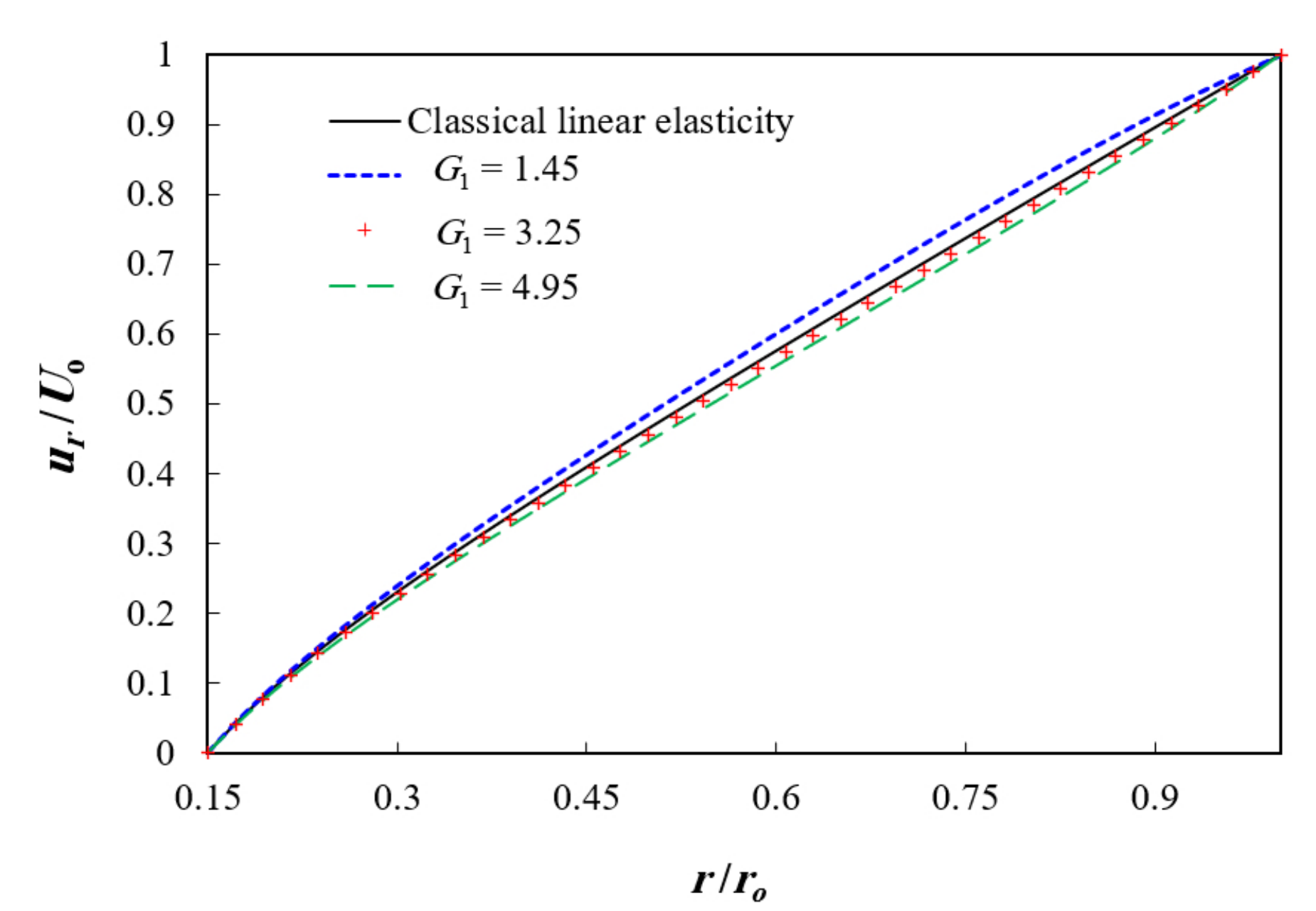}
\caption{Normalized radial displacement of a thick shell with an inner-to-outer radius ratio $\beta=0.15$, comparing the classical elasticity and relaxed micromorphic models for three values of $G_1$ with $r_o/L_{\rm c}=2$, $G_2=5$, $G_3=2$ and $U_i/U_o=0$.}
\label{figure2}
\end{figure}

\begin{figure}[!ht]
\centering
\includegraphics[scale=0.3]{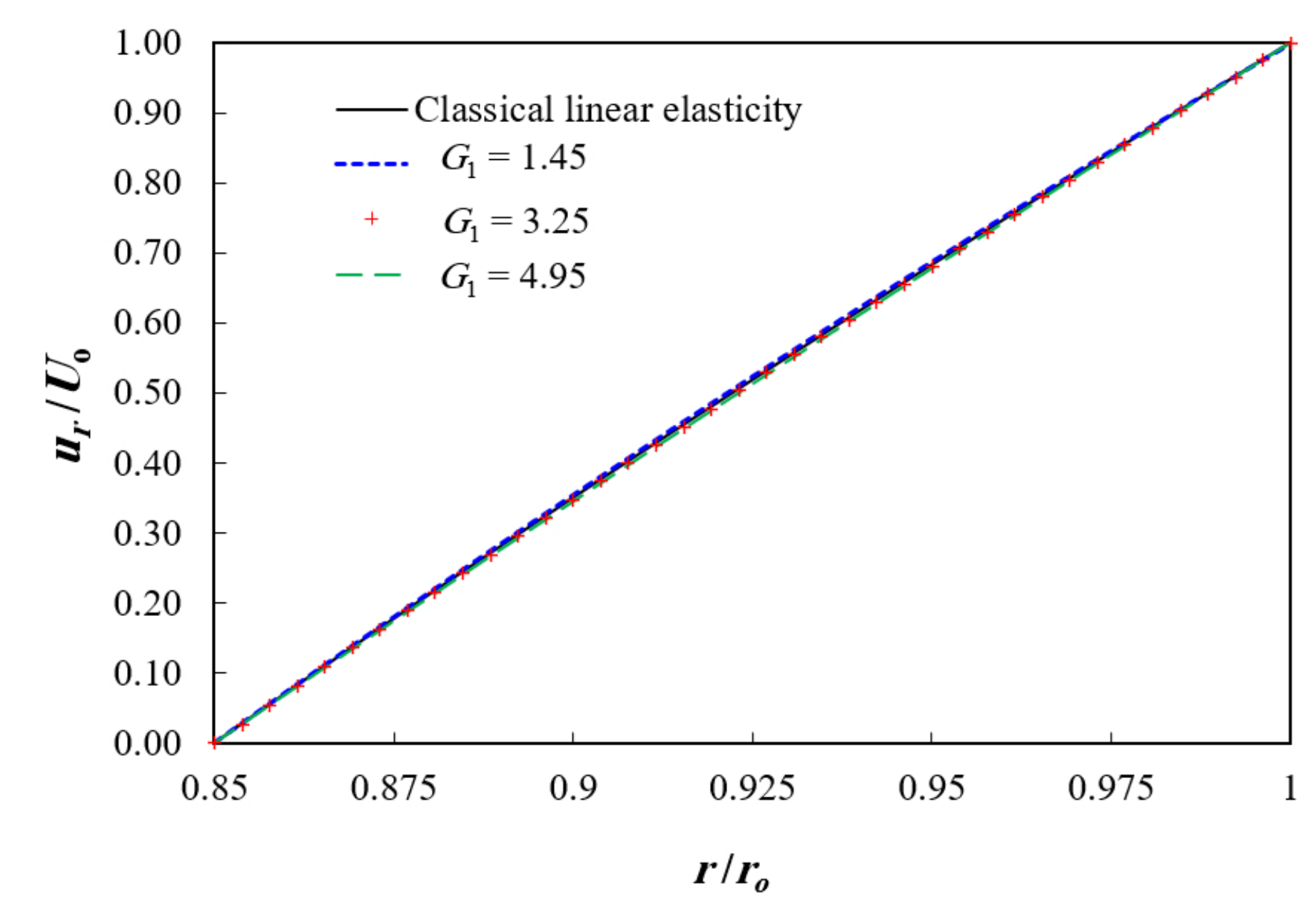}
\caption{Normalized radial displacement of a thin shell with an inner-to-outer radius ratio $\beta=0.85$, comparing the classical elasticity and relaxed micromorphic models for three values of $G_1$ with $r_o/L_{\rm c}=2$, $G_2=5$, $G_3=2$ and $U_i/U_o=0$.}
\label{figure3}
\end{figure}

To examine the effect of $G_2$, numerical results for the thick cylindrical shell are presented in Figure~\ref{figure4}, considering $\beta = 0.15$, $G_1 = 2$, $G_3 = 1.3$, $r_o/L_{\rm c} = 2$, $U_i/U_o = 0$, and $G_2 = 3, 5, 7$. For all considered values of $G_2$, the relaxed micromorphic model predicts larger displacements than classical linear elasticity, with the discrepancy increasing as $G_2$ grows. This behavior is associated with the fact that increasing $G_2$ also causes the parameter $\lambda_e$ to become negative.

\begin{figure}[!ht]
\centering
\includegraphics[scale=0.3]{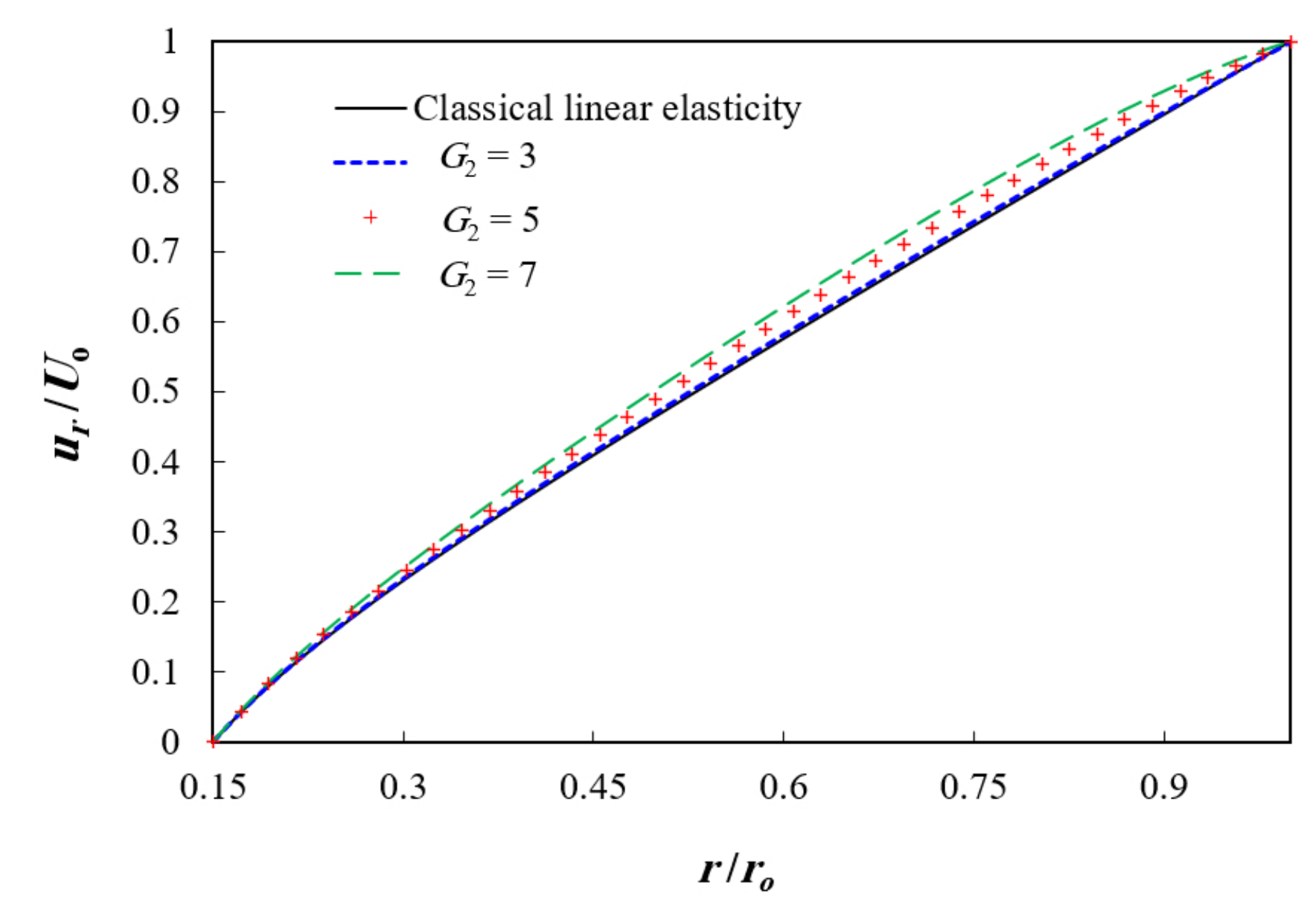}
\caption{Influence of $G_2$ on the normalized radial displacement of a thick cylindrical shell ($\beta = 0.15$) with $G_1 = 2$, $G_3 = 1.3$, $r_o/L_{\rm c} = 2$, and $U_i/U_o = 0$.}
\label{figure4}
\end{figure}

Figure~\ref{figure5} illustrates the normalized radial displacement of the thick cylindrical shell ($\beta=0.15$) for varying values of the micro shear-related parameter $G_3$, with $G_1 = 2.5$, $G_2 = 3.5$, $r_o/L_{\rm c} = 2$, and $U_i/U_o = 0$. For all considered values, the relaxed micromorphic model predicts slightly larger displacements compared to the classical solution, with the magnitude of the deviation increasing as $G_3$ grows. While the overall displacement profile maintains a shape similar to that of classical linear elasticity, the modest upward shift highlights the effect of micromorphic contributions. Collectively, the results presented in Figures~\ref{figure2}, \ref{figure4}, and \ref{figure5} emphasize the significance of accounting for microstructural effects in the mechanical response of thick cylindrical shells.

\begin{figure}[!ht]
\centering
\includegraphics[scale=0.3]{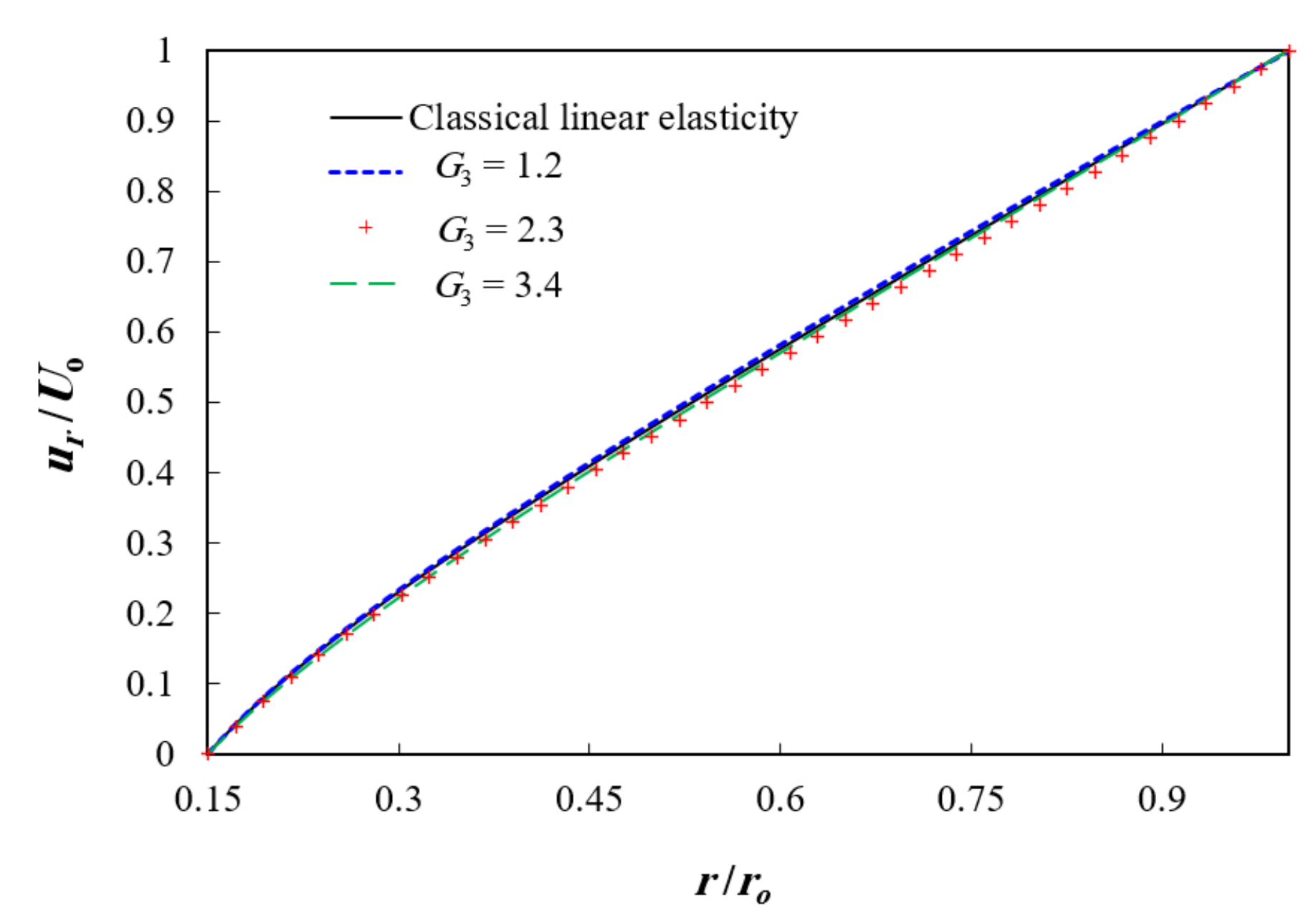}
\caption{Influence of $G_3$ on the normalized radial displacement of a thick cylindrical shell ($\beta = 0.15$) with $G_1 = 2.5$, $G_2 = 3.5$, $r_o/L_{\rm c} = 2$, and $U_i/U_o = 0$.}
\label{figure5}
\end{figure}

Figure~\ref{figure6} illustrates the effect of the inner-to-outer boundary displacement ratio, $\Delta = U_i/U_o$, on the normalized radial displacement of a thick cylindrical shell ($\beta = 0.2$). For these simulations, the material parameters are set to $G_1 = 1.5$, $G_2 = 5.5$, and $G_3 = 2$, with a length scale ratio $r_o/L_c = 1$. The results show that the displacement profile is highly sensitive to the prescribed inner boundary condition. When the inner boundary is moved inward relative to the outer boundary ($\Delta = -0.50$ and $-0.25$), the radial displacement exhibits a steep gradient near the inner surface ($r/r_o = 0.2$). Conversely, as $\Delta$ increases to positive values, the displacement profiles flatten and shift upward. 

\begin{figure}[!ht]
\centering
\includegraphics[scale=0.3]{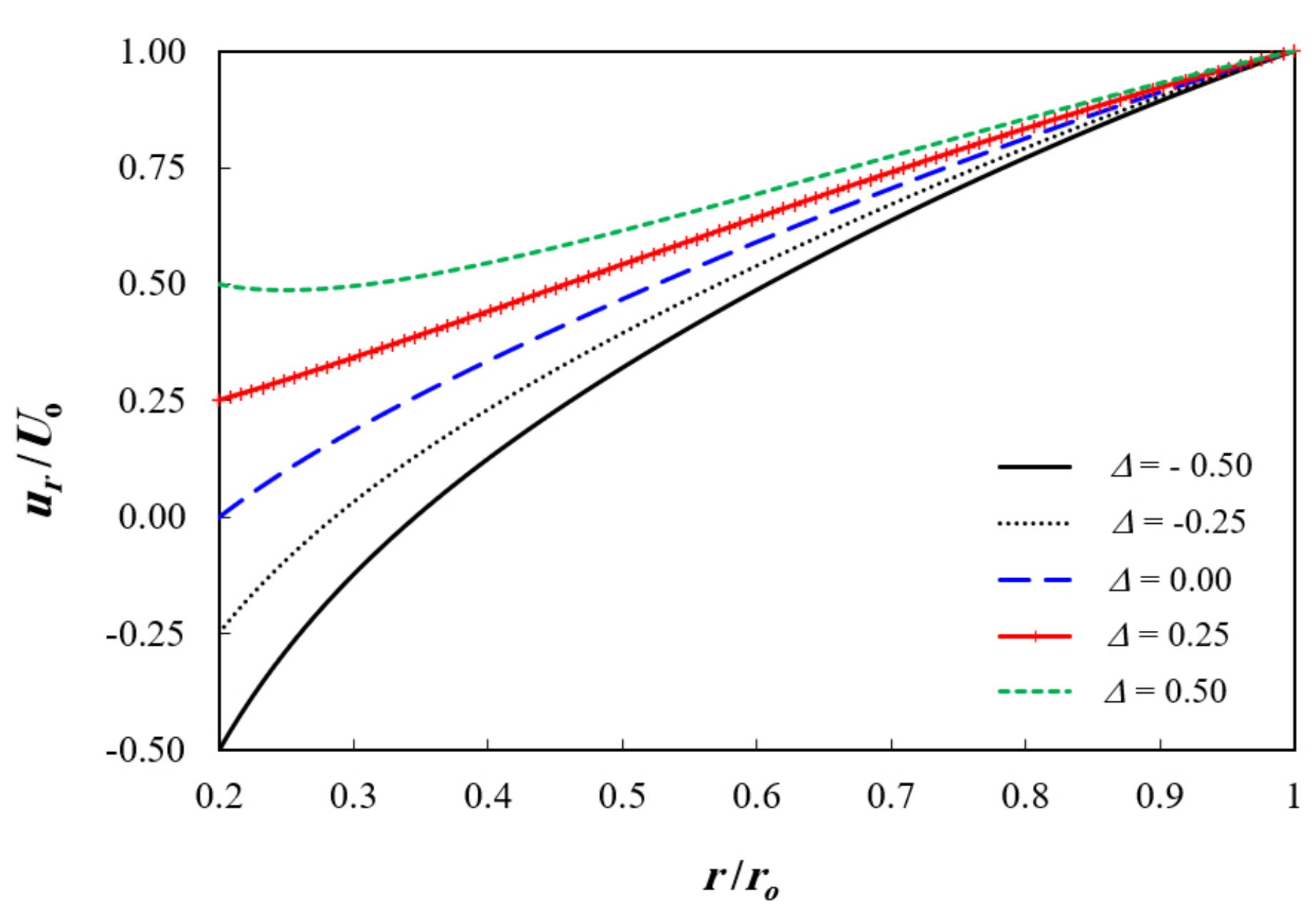}
\caption{Influence of the inner-to-outer boundary displacement ratio $\Delta = U_i/U_o$ on the normalized radial displacement of a thick cylindrical shell ($\beta = 0.2$). The results are plotted for material parameters $G_1 = 1.5$, $G_2 = 5.5$, $G_3 = 2$, and a length scale ratio $r_o/L_c = 1$.}
\label{figure6}
\end{figure}

Finally, to illustrate the effect of the characteristic length, Figures~\ref{figure7} and~\ref{figure8} present the variation of the non-dimensional parameter $\delta$, defined as $[(u_{r}){\text{Micro.}} - (u_{r}){\text{Class.}}]/U_0$, as a function of the radial position. These figures compare the limiting cases $L{c} \rightarrow \infty$ and $L_{c} \rightarrow 0$ for a thick shell with $\beta = 0.25$ and a displacement ratio $\Delta = 0.5$, using the material parameters $G_1 = 1.5$, $G_2 = 5$, and $G_3 = 2$. The results clearly show that the characteristic length strongly influences the displacement response. Figure~\ref{figure7} indicates that as the ratio $r_o/L_c$ decreases (corresponding to the limit $L_c \rightarrow \infty$), the deviation $\delta$ approaches zero throughout the entire radial domain. Likewise, Figure~\ref{figure8} demonstrates that for very large values of $r_o/L_c$ (corresponding to the limit $L_c \rightarrow 0$), the micromorphic solution again converges to the classical elasticity solution. Therefore, in both limiting cases of very large and very small characteristic lengths, the relaxed micromorphic model effectively reduces to the classical model, indicating a negligible influence of microstructural effects.

\begin{figure}[!ht]
\centering
\includegraphics[scale=0.3]{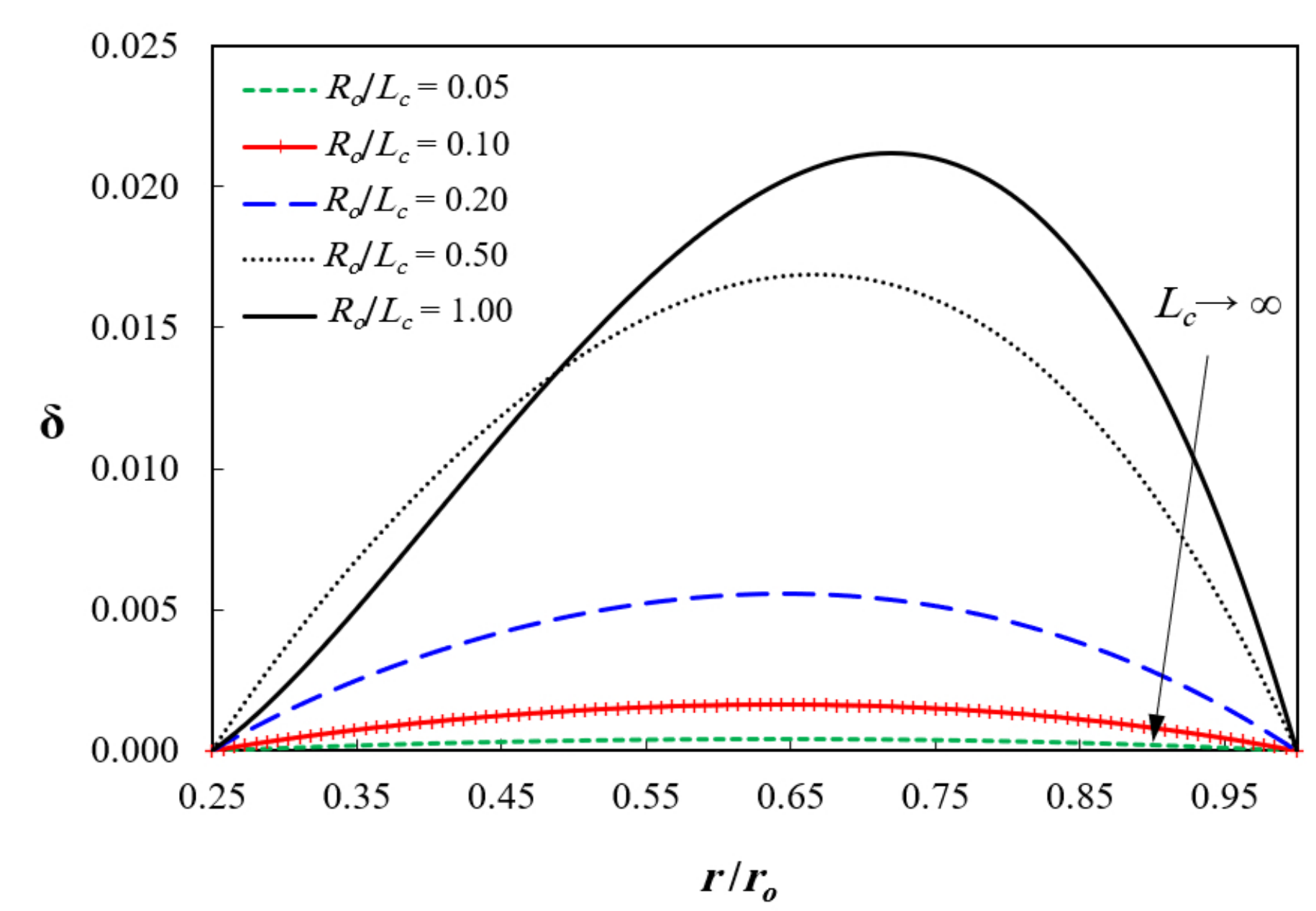}
\caption{Comparison between the relaxed micromorphic and the classical models for $R/L_{c}=\{{0.05,0.10,0.20,0.50,1.00}\}$.}
\label{figure7}
\end{figure}

\begin{figure}[!ht]
\centering
\includegraphics[scale=0.3]{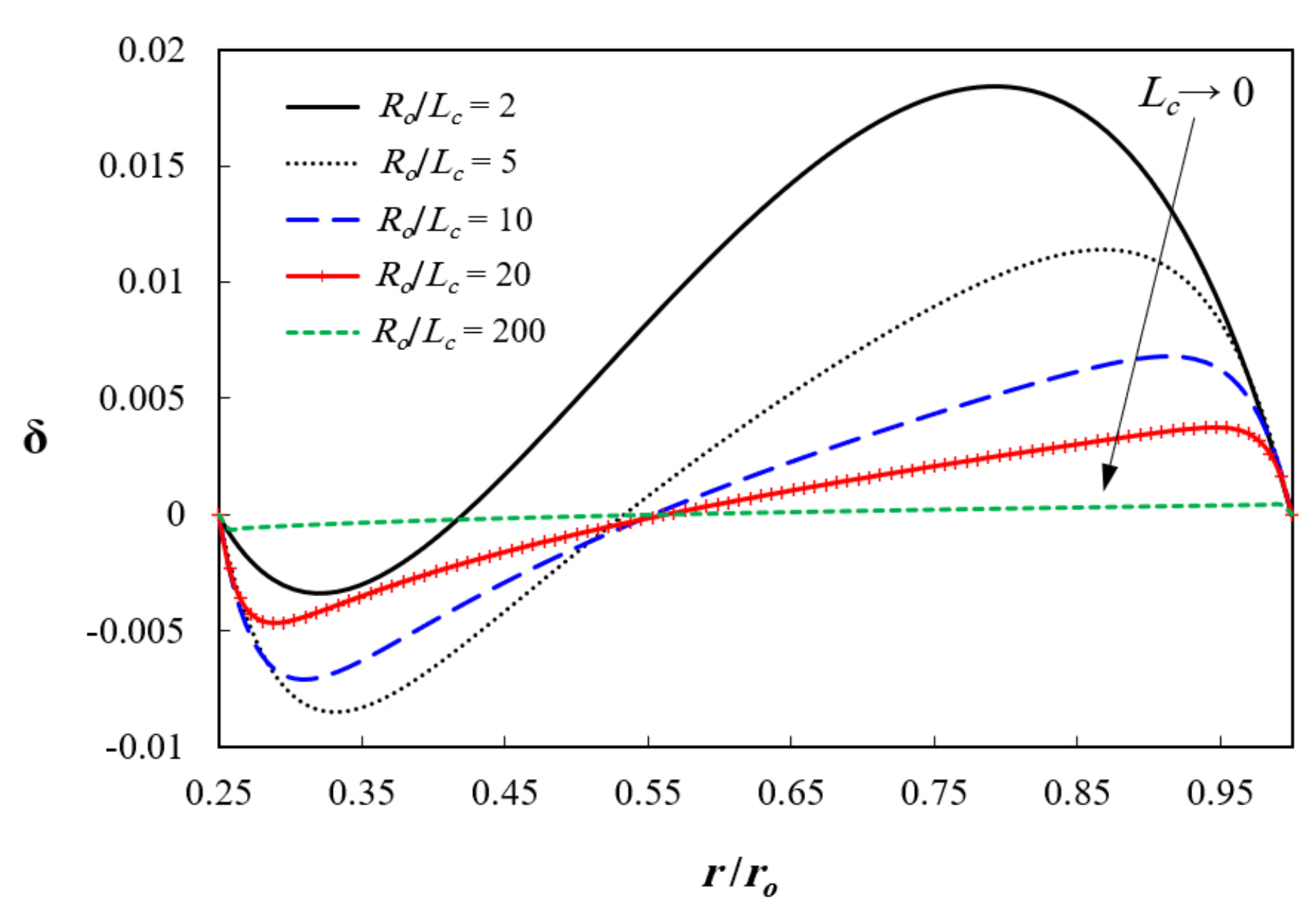}
\caption{Comparison between the relaxed micromorphic and the classical models for $R/L_{c}=\{{2,5,10,20,200}\}$.}
\label{figure8}
\end{figure}

\section{Conclusions}
This work has developed an exact analytical solution for the elastostatic deformation of long cylindrical shells within the framework of the isotropic relaxed micromorphic continuum. By exploiting the axisymmetric nature of the problem and introducing appropriate auxiliary variables, the governing equations were reduced to a modified Bessel differential equation. Closed-form expressions were obtained for the radial displacement field as well as the relevant components of the micro-distortion tensor. The analytical solution reveals the fundamental role of the characteristic length parameter in governing the mechanical response of microstructured materials. In contrast to classical elasticity, the relaxed micromorphic model predicts size-dependent behavior arising from microstructural interactions. The solution naturally recovers the classical elasticity response in the limiting case of vanishing characteristic length, thereby confirming the consistency of the formulation.

The results establish a rigorous benchmark solution that can be used to verify finite element implementations and numerical algorithms for relaxed micromorphic continua. The present work contributes to the theoretical understanding of microstructured cylindrical shells. It supports the broader application of relaxed micromorphic models in the analysis and design of advanced materials and structures exhibiting size-dependent behavior. Future work may extend the present formulation to dynamic loading conditions, layered or functionally graded materials, and more general boundary conditions, thereby further expanding the applicability of relaxed micromorphic theory to engineering problems involving microstructured solids.%

\begingroup
\setstretch{1}
\setlength\bibitemsep{3pt}
\printbibliography
\endgroup

\end{document}